\documentclass[12pt]{article}
\usepackage{amsmath,amssymb}

\newcommand{\1}{^{-1}}
\newcommand{\tr}{{}^t\!}

\newcommand{\iso}{\cong}
\newcommand{\ot}{\leftarrow}
\newcommand{\into}{\hookrightarrow}
\newcommand{\lodot}{_\bullet}
\newcommand{\tensor}{\otimes}
\newcommand{\J}{\left(\begin{smallmatrix} 0&I \\ I&0\end{smallmatrix}\right)}
\newcommand{\dgn}{\mathrm{dgn}}
\newcommand{\aff}{\mathbb A}
\newcommand{\kk}{\mathbf k}

\newcommand{\sP}{\mathcal P}
\newcommand{\Oh}{\mathcal O}
\newcommand{\sS}{\mathcal S}
\newcommand{\C}{\mathbb C}
\newcommand{\Gm}{{\mathbb G}_m}
\newcommand{\PP}{\mathbb P}
\newcommand{\Z}{\mathbb Z}

\newcommand{\al}{\alpha}
\newcommand{\Ga}{\Gamma}
\newcommand{\De}{\Delta}
\newcommand{\ka}{\kappa}
\newcommand{\la}{\lambda}
\newcommand{\fie}{\varphi}
\newcommand{\si}{\sigma}
\newcommand{\Si}{\Sigma}
\newcommand{\om}{\omega}

\newcommand{\Span}[1]{\left<#1\right>}
\DeclareMathOperator{\Adj}{Adj}
\DeclareMathOperator{\coker}{coker}
\DeclareMathOperator{\codim}{codim}
\DeclareMathOperator{\depth}{depth}

\DeclareMathOperator{\wt}{wt}
\DeclareMathOperator{\Ext}{Ext}
\DeclareMathOperator{\GL}{GL}
\DeclareMathOperator{\Gr}{Gr}
\DeclareMathOperator{\Hom}{Hom}
\DeclareMathOperator{\Frac}{Frac}

\DeclareMathOperator{\Mat}{Mat}
\DeclareMathOperator{\Mor}{Mor}
\DeclareMathOperator{\Pf}{Pf}
\DeclareMathOperator{\OGr}{OGr}

\DeclareMathOperator{\Orth}{O}
\DeclareMathOperator{\Pic}{Pic}
\DeclareMathOperator{\Pin}{Pin}
\DeclareMathOperator{\Segre}{Segre}
\DeclareMathOperator{\Spec}{Spec}
\DeclareMathOperator{\SpH}{SpH}
\DeclareMathOperator{\Spin}{Spin}
\DeclareMathOperator{\Supp}{Supp}
\DeclareMathOperator{\Sym}{Sym}
\newtheorem{theorem}{Theorem}[section]
\newtheorem{cla}[theorem]{Claim}
\newtheorem{exa}[theorem]{Example}
\newtheorem{exc}[theorem]{Exercise}
\newtheorem{lem}[theorem]{Lemma}
\newtheorem{conj}[theorem]{Conjecture}

\newtheorem{propdef}[theorem]{Proposition--Definition}
\numberwithin{equation}{section}
\newcommand{\QED}{\ifhmode\unskip\nobreak\fi\quad\ensuremath{\mathrm{QED}}}

\title{Gorenstein in codimension 4 -- \\ the general structure theory}
\author{Miles Reid}
\date{}

\begin{document}
\maketitle

\begin{abstract}

I describe the projective resolution of a codimension~4 Gorenstein ideal,
aiming to extend Buchsbaum and Eisenbud's famous result in
codimension~3. The main result is a structure theorem stating that the
ideal is determined by its $(k+1)\times2k$ matrix of first syzygies,
viewed as a morphism from the ambient regular space to the Spin-Hom
variety $\SpH_k\subset\Mat(k+1,2k)$. This is a general result
encapsulating some theoretical aspects of the problem, but, as it
stands, is still some way from tractable applications.

\end{abstract}

This paper introduces the {\em Spin-Hom varieties}
$\SpH_k\subset\Mat(k+1,2k)$ for \hbox{$k\ge3$}, that I define as almost
homo\-geneous spaces under the group $\GL(k+1)\times\Orth(2k)$ (see
\ref{s!SpH}). These serve as key varieties for the $(k+1)\times 2k$ first
syzygy matrixes of codimension~4 Gorenstein ideals $I$ in a polynomial
ring $S$ plus appropriate presentation data; the correspondence takes
$I$ to its matrix of first syzygies. Such ideals $I$ are parametrised by
an open subscheme of $\SpH_k(S)=\Mor(\Spec S,\SpH_k)$. The open
condition comes from the Buchsbaum--Eisenbud exactness criterion ``What
makes a complex exact?'' \cite{BE1}: the classifying map $\al\colon\Spec
S\to\SpH_k$ must hit the {\em degeneracy locus} of $\SpH_k$ in
codimension $\ge4$.

The map $\al$ has {\em Cramer-spinor coordinates} $L_i$ and $\si_J$ in
standard representations $\kk^{k+1}$ and $\kk^{2^{k-1}}$ of $\GL(k+1)$
and $\Pin(2k)$ (see \ref{s!sc}), and the $k\times k$ minors of $M_1(I)$
are in the product ideal $I\cdot\Sym^2(\{\si_J\})$. The spinors
themselves should also be in $I$, so that the $k\times k$ minors of
$M_1(I)$ are in $I^3$; this goes some way towards explaining the
mechanism that makes the syzygy matrix $M_1(I)$ ``drop rank by~3 at one
go'' -- it has rank $k$ outside $V(I)=\Spec(S/I)$ and $\le k-3$ on
$V(I)$.

The results here are not yet applicable in any satisfactory way, and
raise almost as many questions as they answer. While Gorenstein
codimension~4 ideals are subject to a structure theorem, that I believe
to be the correct codimension~4 generalisation of the famous
Buchsbaum--Eisenbud theorem in codimension~3 \cite{BE2}, I do not say
that this makes them tractable.

\medskip\noindent{\bf Thanks } I am grateful to Chen Jungkai for
inviting me to AGEA, and to him and his colleagues at University of
Taiwan for generous hospitality. My visit was funded by Korean
Government WCU Grant R33-2008-000-10101-0, which also partly funded my
work over the last 4 years, and I am very grateful to Lee Yongnam for
setting up the grant and administering every aspect of it. I wish to
thank Fabrizio Catanese, Eduardo Dias, Sasha Kuznetsov and Liu Wenfei
for contributing corrections, questions and stimulating discussion. I
owe a particular debt of gratitude to Alessio Corti for detailed
suggestions that have helped me improve the layout and contents of the
paper.

\medskip\noindent{\bf Website } See
www.warwick.ac.uk/staff/Miles.Reid/codim4 for material
accompanying this paper.

\section{Introduction} \label{s!1}

Gorenstein rings are important, appearing throughout algebra, algebraic
geometry and singularity theory. A common source is Zariski's standard
construction of graded ring over a polarised variety $X,L$: the graded
ring $R(X,L)=\bigoplus_{n\ge0}H^0(X,nL)$ is a Gorenstein ring under
natural and fairly mild conditions (cohomology vanishing plus $K_X=k_XL$
for some $k_X\in\Z$, see for example \cite{GW}). Knowing how to
construct $R(X,L)$ by generators and relations gives precise answer to
questions on embedding $X\into\PP^n$ and determining the equations of
the image.

\subsection{Background and the Buchsbaum--Eisenbud result}
I work over a field $\kk$ containing $\frac12$ (such as $\kk=\C$, but
see \ref{ss!J} for the more general case). Let $S=\kk[x_1,\dots,x_n]$ be
a positively graded polynomial ring with $\wt x_i=a_i$, and $R=S/I_R$ a
quotient of $S$ that is a Gorenstein ring. Equivalently, $\Spec
R\subset\Spec S=\aff^n_{\kk}$ is a Gorenstein graded scheme. By the
Auslander--Buchsbaum form of the Hilbert syzygies theorem, $R$ has a
minimal free graded resolution $P\lodot$ of the form
\begin{equation}
\begin{aligned}
0 \ot{} &P_0 \ot P_1 \ot\cdots \ot P_c\ot 0 \\
&\downarrow \\
&R
\end{aligned}
\label{eq!resP}
\end{equation}
where $P_0=S\to R=S/I_R$ is the quotient map, and $P_1\to S$ gives a
minimum set of generators of the ideal $I_R$. Here the length $c$ of the
resolution equals $n-\depth R$, and each $P_i$ is a graded free module
of rank $b_i$. I write $P_i=b_iS$ (as an abbreviation for $S^{\oplus
b_i})$, or $P_i=\bigoplus_{j=1}^{b_i} S(-d_{ij})$ if I need to keep
track of the gradings. The condition $\depth R=\dim R$ that the depth is
maximal characterises the Cohen--Macaulay case, and then $c=\codim
R=\codim(\Spec R\subset\Spec S)$. If in addition $P_c$ is a free module
of rank~1, so that $P_c\iso S(-\al)$ with $\al$ the {\em adjunction
number}, then $R$ is a Gorenstein ring of canonical weight
$\ka_R=\al-\sum a_i$; for my purposes, one can take this to be the
definition of Gorenstein.

Duality makes the resolution \eqref{eq!resP} symmetric: the dual complex
$(P\lodot)^\vee = \Hom_S(P\lodot,P_c)$ resolves the dualising module
$\om_R=\Ext^c_S(R,\om_S)$, which is isomorphic to $R$ (or, as a graded
module, to $R(\ka_R)$ with $\ka_R=\al-\sum a_i$), so that $P\lodot\iso
(P\lodot)^\vee$. In particular the Betti numbers $b_i$ satisfy the
symmetry $b_{c-i}=b_i$, or
\[
P_{c-i}=\Hom_S(P_i,P_c) \iso\bigoplus_{j=1}^{b_i} S(-\al+d_{ij}),
\quad\hbox{where}\quad P_i=\bigoplus_{j=1}^{b_i} S(-d_{ij}).
\]

The Buchsbaum--Eisenbud symmetriser trick \cite{BE2} adds precision to
this (this is where the assumption $\frac12\in S$ comes into play):
\begin{quote}
{\em There is a symmetric perfect pairing $S^2(P\lodot)\to P_c$ inducing
the duality $P\lodot\iso (P\lodot)^\vee$.}
\end{quote}
The idea is to pass from $P\lodot$ as a resolution of $R$ to the complex
$P\lodot\tensor P\lodot$ (the total complex of the double complex) as a
resolution of $R\tensor_S R$ (left derived tensor product), then to
replace $P\lodot\tensor P\lodot$ by its symmetrised version
$S^2(P\lodot)$. In the double complex $P\lodot\tensor P\lodot$, one
decorates the arrows by signs $\pm1$ to make each rectangle anticommute
(to get $d^2=0$). The symmetrised complex $S^2(P\lodot)$ then involves
replacing the arrows by half the sum or differences of symmetrically
placed arrows. (This provides lots of opportunities for confusion about
signs!)

For details, see \cite{BE2}. The conclusion is that $P\lodot$ has a
$\pm$-symmetric bilinear form that induces perfect pairings $P_i\tensor
P_{c-i}\to P_c=S$ for each $i$, compatible with the differentials.

The Buchsbaum--Eisenbud structure theorem in codimension~3 is a simple
consequence of this symmetry, and a model for what I try to do in this
paper. Namely, in codimension~3 we have
\begin{equation}
0\ot P_0 \ot P_1 \ot P_2 \ot P_3 \ot 0,
\label{eq!cod3}
\end{equation}
with $P_0=S$, $P_3\iso S$, $P_2=\Hom(P_1,P_3)\iso P_1^\vee$, and the
matrix $M$ defining the map $P_1\ot P_2$ is skew (that is,
antisymmetric). If I set $P_1=nS$ then the respective ranks of the
differentials in \eqref{eq!cod3} are $1$, $n-1$ and $1$; since $M$ is
skew, his rank must be even, so that $n=2\nu+1$. Moreover, the kernel
and cokernel are given by the Pfaffians of $M$, by the skew version of
Cramer's rule.\goodbreak

Generalising the Buchsbaum--Eisenbud Theorem to codimension~4 has been a
notoriously elusive problem since the 1970s.

\subsection{Main aim} \label{s!ma}
This paper starts by describing the shape of the resolution of a
codimension~4 Gorenstein ring by analogy with \eqref{eq!cod3}. The first
syzygy matrix $M_1\colon P_1\ot P_2$ is a $(k+1)\times2k$ matrix whose
$k+1$ rows generically span a maximal isotropic space of the symmetric
quadratic form on $P_2$. The ideal $I_R$ is generated by the entries of
the map $L\colon P_0\ot P_1$, which is determined by the linear algebra
of quadratic forms as the linear relation that must hold between the
$k+1$ rows of $M_1$.

This is all uncomplicated stuff, deduced directly from the symmetry
trick of \cite{BE2}. It leads to the definition of the Spin-Hom
varieties $\SpH_k$ in the space of $(k+1)\times2k$ matrixes (see
Section~\ref{s!SpH}). The first syzygy matrix $M_1$ is then an
$S$-valued point of $\SpH_k$, or a morphism $\al\colon\Spec S\to\SpH_k$.

The converse is more subtle, and is the main point of the paper. By
construction, $\SpH_k$ supports a short complex $\sP_1\ot\sP_2\ot\sP_3$
of free modules with a certain universal property. If we were allowed to
restrict to a smooth open subscheme $S^0$ of $\SpH_k$ meeting the
degeneracy locus $\SpH^{\dgn}_k$ in codimension~4, the reflexive hull of
the cokernel of $M_1$ and the kernel of $M_2$ would provide a complex
$\sP\lodot$ that resolves a sheaf of Gorenstein codimension~4 ideals in
$S^0$. (This follows by the main proof below).

Unfortunately, this is only an adequate description of codimension~4
Gorenstein ideals in the uninteresting case of complete intersection
ideals. Any other case necessarily involves smaller strata of $\SpH_k$,
where $\SpH_k$ is singular. Thus to cover every codimension~4 Gorenstein
ring, I am forced into the logically subtle situation of a universal
construction whose universal space does not itself support the type of
object I am trying to classify, namely a Gorenstein codimension~4 ideal.
See \ref{s!sbt} for further discussion of this point.

Main Theorem~\ref{s!main} gives the universal construction. To
paraphrase: for a polynomial ring $S$ graded in positive degrees, there
is a 1-to-1 correspondence between:
\begin{enumerate}
\renewcommand{\labelenumi}{(\arabic{enumi})}

\item Gorenstein codimension 4 graded ideals $I\subset S$ and

\item graded morphisms $\al\colon\Spec S\to\SpH_k$ for which
$\al\1(\SpH^{\dgn}_k)$ has codimension $\ge4$ in $\Spec S$.

\end{enumerate}
I should say at once that this is intended as a theoretical structure
result. It has the glaring weakness that it does not so far make any
tractable predictions even in model cases (see \ref{s!app} for a
discussion). But it is possibly better than no structure result at all.

\subsection{Contents of the paper} \label{s!cont}

Section~\ref{s!fr} describes the shape of the free resolution and its
symmetry, following the above introductory discussion.
Section~\ref{s!SpH} defines the Spin-Hom variety
$\SpH_k\subset\Mat(k+1,2k)$, to serve as my universal space. The
definition takes the form of a quasihomogeneous space for the complex
Lie group $G=\GL(k+1)\times\Orth(2k)$ or its spin double cover
$\GL(k+1)\times \Pin(2k)$. More explicitly, define $\SpH_k$ as the
closure of the $G$-orbit $\SpH_k^0=G\cdot M_0$ of the {\em typical
matrix} $M_0=(\begin{smallmatrix}I_k&0\\0&0\end{smallmatrix})$ under the
given action of $G=\GL(k+1)\times\Orth(2k)$ on \hbox{$\Mat(k+1,2k)$}.

The degeneracy locus $\SpH^{\dgn}_k$ is the complement
$\SpH_k{}\setminus{}\SpH_k^0$. Once these definitions are in place,
Section~\ref{s!main} states the main theorem, and proves it based on the
exactness criterion of \cite{BE1}.

The Spin-Hom varieties $\SpH_k$ have a rich structure arising from
representation theory. A matrix $M_1\in\SpH_k^0$ can be viewed as an
isomorphism between a $k$-dimensional space in $\kk^{k+1}$ and a maximal
isotropic space for $\fie$ in $\kk^{2k}$. This displays $\SpH_k^0$ as a
principal $\GL(k)$ bundle over $\PP^k\times\OGr(k,2k)$.
Section~\ref{s!3} discusses the properties of the $\SpH_k$ in more
detail, notably their symmetry under the maximal torus and Weyl group.
The spinor and nonspinor sets correspond to the two different spinor
components $\OGr(k,2k)$ and $\OGr'(k,2k)$ of the maximal isotropic
Grassmannian.

I introduce the Cramer-spinor coordinates $\si_J$ in \ref{s!sc}; the
main point is that, for a spinor subset $J\cup J^c$, the $(k+1)\times k$
submatrix of $M_1\in\SpH_k$ formed by those columns has top wedge
factoring as $(L_1,\dots,L_{k+1})\cdot\si_J^2$ where $L\colon P_0\ot
P_1$ is the vector of equations (see Lemma~\ref{s!spin}). Ensuring
that the appropriate square root $\si_J$ is defined as an element
$\si_J\in S$ involves the point that, whereas the spinor bundle defines
a 2-torsion Weil divisor on the affine orthogonal Grassmannian
$a\OGr(k,2k)\subset\bigwedge^k\kk^{2k}$ (the affine cone over
$\OGr(k,2k)$ in Pl\"ucker space) and on $\SpH_k$, its birational
transform under the classifying maps $\al\colon\Spec S\to\SpH_k$ of
Theorem~\ref{s!main} is the trivial bundle on $\Spec S$.

The spinor coordinates vanish on the degeneracy locus $\SpH^{\dgn}_k$
and define an equivariant morphism
$\SpH_k^0\to\kk^{k+1}\tensor\kk^{2^{k-1}}$. At the same time, they
vanish on the nonspin variety $\SpH'_k$, corresponding to the other
component $\OGr'(k,2k)$ of the Grassmannian of maximal isotropic
subspaces; this has nonspinor coordinates, that vanish on $\SpH_k$.
Between them, these give set theoretic equations for $\SpH_k$ and its
degeneracy locus.\goodbreak

The final Section~\ref{s!4} discusses a number of issues with my
construction and some open problems and challenges for the future.

\section{The main result} \label{s!2}

For a codimension~4 Gorenstein ideal $I$ with $k+1$ generators, the
module $P_2$ of first syzygies is a $2k$ dimensional orthogonal space
with a nondegenerate (symmetric) quadratic form $\fie$. The $k+1$ rows
of the first syzygy matrix $M_1(R)$ span an isotropic subspace in $P_2$
with respect to $\fie$. Since the maximal isotropic subspaces are
$k$-dimensional, this implies a linear dependence relation
$(L_1,\dots,L_{k+1})$ that bases $\coker M_1$ and thus provides the
generators of $I$. A first draft of this idea was sketched in \cite{Ki},
10.2.

\subsection{The free resolution} \label{s!fr}
Let $S=\kk[x_1,\dots,x_N]$ be the polynomial ring over an algebraically
closed field $\kk$ of characteristic $\ne2$, graded in positive degrees.
Let $I_R$ be a homogeneous ideal with quotient $R=S/I_R$ that is
Gorenstein of codimension~4; equivalently, $I_R$ defines a codimension~4
Gorenstein graded subscheme
\[
V(I_R)=\Spec R\subset\aff^N_\kk=\Spec S.
\]
Suppose that $I_R$ has $k+1$ generators $L_1,\dots,L_{k+1}$. It follows
from the Auslander--Buchsbaum form of the Hilbert syzygies theorem and
the symmetriser trick of Buchsbaum--Eisenbud \cite{BE2} that the free
resolution of $R$ is
\begin{equation}
0 \ot P_0 \ot P_1 \ot P_2 \ot P_3 \ot P_4 \ot 0,
\label{eq!P4}
\end{equation}
where $P_0=S$, $P_4\iso S$, $P_3=\Hom(P_1,P_4)\iso P_1^\vee$; and
moreover, $P_2$ has a nondegenerate {\em symmetric} bilinear form
$\fie\colon S^2P_2\to P_4$ compatible with the complex $P\lodot$, so
that $P_2\to P_1$ is dual to $P_3\to P_2$ under $\fie$. The simple cases
of \ref{s!exa}, Examples~\ref{exa!Kos4}--\ref{exa!T} give a sanity
check (just in case you are sceptical about the symmetry of $\fie$).

A choice of basis of $P_2$ gives $\fie$ the standard block
form\footnote{In the graded case this is trivial because $\fie$ is
homogeneous of degree 0, so is basically a nondegenerate quadratic form
on a vector space $V_2$ with $P_2=V_2\tensor S$. See the discussion in
\ref{ss!J} for the more general case.} $\J$. Then the first syzygy matrix
in \eqref{eq!P4} is $M_1(R)=(A\,B)$, where the two blocks are
$(k+1)\times k$ matrixes satisfying
\begin{equation}
(A\,B) \J \,^t(A\,B)=0,
\label{eq!AtB}
\end{equation}
that is, $A \,\tr B + B \,\tr A = 0$, or $A \,\tr B$ is {\em skew}.
I call this a $(k+1)\times2k$ resolution (meaning that the defining
ideal $I_R$ has $k+1$ generators yoked by $2k$ first syzygies).

The number of equations in \eqref{eq!AtB} is $\binom {k+2}2$. For
example, in the typical case $k=8$, the variety defined by
\eqref{eq!AtB} involves $\binom {k+2}2=45$ quadratic equations in
$2k(k+1)=144$ variables. The scheme $V_k$ defined by \eqref{eq!AtB}
appears in the literature as the {\em variety of complexes}. However it
is not really the right object -- it breaks into 2 irreducible
components for spinor reasons, and it is better to study just one, which
is my $\SpH_k$.

\subsection{The general fibre} \label{s!gf}
Let $\xi\in\Spec S=\aff^N$ be a point outside $V(I_R)=\Spec R$ with
residue field $K=\kk(\xi)$ (for example, a $\kk$-valued point, with $K=\kk$,
or the generic point, with $K=\Frac S$). Evaluating \eqref{eq!P4} at
$\xi$ gives the exact sequence of vector spaces
\begin{equation}
0 \ot V_0 \ot V_1 \ot V_2 \ot V_3 \ot V_4 \ot 0
\label{eq!Vi}
\end{equation}
over $K$, where $V_0=K$, $V_4\iso K$, $V_1=(k+1)K$, $V_3=\Hom(V_1,V_4)\iso
V_1^\vee$, and $V_2=2kK$ with the nondegenerate quadratic form
$\fie=\J$. Over $K$, the maps in \eqref{eq!Vi} can be written as the
matrixes
\begin{equation}
\begin{pmatrix} 0 & \dots & 0 & 1
\end{pmatrix}
\begin{pmatrix} I_k & 0 \\ 0 & 0
\end{pmatrix}
\begin{pmatrix} 0 & 0 \\ I_k & 0
\end{pmatrix}
\left(
\begin{smallmatrix} 0 \\[-3pt] \vdots \\[3pt] 0 \\ 1
\end{smallmatrix} \right).
\label{eq!genf}
\end{equation}
This data determines a fibre bundle over $\aff^N\setminus V(I_R)$ with
the exact complex \eqref{eq!Vi} as fibre, and structure group the
orthogonal group of the complex, which I take to be
$\GL(k+1)\times\Orth(2k)$ or its double cover $\GL(k+1)\times\Pin(2k)$.

\subsection{Simple examples} \label{s!exa}
\begin{exa}\label{exa!Kos4}\rm A codimension 4 complete intersection
has $L=(x_1,x_2,x_3,x_4)$ and Koszul syzygy matrix
\begin{equation}
(A\,B) =
\begin{pmatrix} -x_4 & . & . && . & x_3 & -x_2 \\
. & -x_4 & . && -x_3 & . & x_1 \\
. & . & -x_4 && x_2 & -x_1 & . \\
x_1 & x_2 & x_3 && . & . & .
\end{pmatrix}.
\label{eq!Kos4}
\end{equation}
In this choice, $A=M_{1,2,3}$ has rank 3 and
$\bigwedge^3A=x_4^2\cdot(x_1,\dots,x_4)$. See \ref{s!sc} for spinors. A
spinor subset $J\cup J^c$ has an odd number $i$ of columns from $A$ and
the complementary $3-i$ columns from $B$. For example, columns $1,5,6$
give a $4\times3$ matrix with $\bigwedge^3
M_{1,5,6}=x_1^2\cdot(x_1,x_2,x_3,x_4)$.
\end{exa}

\begin{exa}\label{exa!PfHy}\rm Another easy case is that of a
hypersurface section $h=0$ in a codimension~3 ideal given by the
Pfaffians $\Pf_i$ of a $(2l+1)\times(2l+1)$ skew matrix $M$. The
syzygy matrix is
\begin{equation}
(A\,B) =
\begin{pmatrix} -hI_{2l+1} && M \\[6pt]
\Pf_1 \dots \Pf_{2l+1} && 0\dots0
\end{pmatrix}.
\label{eq!Pfh}
\end{equation}
One sees that a spinor $\si_J$ corresponding to $2l+1-2i$ columns from
$A$ and a complementary $2i$ from $B$ is of the form $h^{l-i}$ times a
diagonal $2i\times2i$ Pfaffian of $M$. Thus the top wedge of the left-hand
block $A$ of \eqref{eq!Pfh} equals $\si^2\cdot(h,\Pf_1,\dots,\Pf_{2l+1})$
where $\si=h^l$.

\end{exa}

\begin{exa}\label{exa!T} \rm The extrasymmetric matrix
\begin{equation}
M = \begin{pmatrix}
a&b&&d&e&f \\
&c&&e&g&h \\
&&& f&h&i \\[6pt]
&&&& -\la a& -\la b \\
&&&&& -\la c
\end{pmatrix}
\end{equation}
with a single multiplier $\la$ is the simplest case of a Tom
unprojection (see \cite{TJ}, Section~9 for details). Let $I$ be the
ideal generated by the $4\times4$ Pfaffians of $M$. The diagonal entries
$d,g,i$ of the $3\times3$ symmetric top right block are all unprojection
variables; thus $i$ appears linearly in 4 equations of the form
$i\cdot(a,d,e,g)=\cdots$, and eliminating it projects to the
codimension~3 Gorenstein ring defined by the Pfaffians of the top left
$5\times5$ block.

If $\la\in S$ is a perfect square, $I$ is the ideal of
$\Segre(\PP^2\times\PP^2)\subset\PP^8$ up to a coordinate change, but
the Galois symmetry $\sqrt\la\mapsto-\sqrt\la$ swaps the two factors.
See \cite{TJ}, Section~9 for more details, and for several more families
of examples; in any of these cases, writing out the resolution matrixes
$(A\,B)$ with the stated isotropy property makes a demanding but
rewarding exercise for the dedicated student.

By extrasymmetry, out of the 15 entries of $M$, 9 are independent and 6
repeats. His $4\times4$ Pfaffians follow a similar pattern. I write the
9 generators of the ideal $I$ of Pfaffians as the vector $L=$
\begin{equation*}
\begin{array}{l}
\bigl[ \la ac+eh-fg, \enspace -\la ab-dh+ef, \enspace \la a^2+dg-e^2, \\
\kern5em ah-bg+ce, \enspace -af+be-cd, \enspace \la b^2+di-f^2, \\
\kern10em \la bc+ei-fh, \enspace \la c^2+gi-h^2, \enspace ai-bh+cf \bigr]
\end{array}
\end{equation*}

Its matrix of first syzygies $M_1$ is the transpose of
\begin{equation}
\begin{matrix}
. & a & b & d & e & . & . & . & . \\
-a & . & c & e & g & . & . & . & . \\
-b & -c & . & f & h & . & . & . & . \\
-d & -e & -f & . & -\la a & . & . & . & . \\
-e & -g & -h & \la a & . & . & . & . & . \\[6pt]
-h & . & . & \la c & . & . & g & -e & . \\
f & -h & . & -\la b & \la c & -g & . & d & . \\
. & f & . & . & -\la b & e & -d & . & . \\[12pt]
i & . & . & . & . & . & -h & f & -\la c \\
. & i & . & . & . & h & -f & . & \la b \\
. & h & i & . & -\la c & . & e & -d & -\la a \\
. & . & . & i & . & . & -c & b & -h \\
. & . & . & . & i & c & -b & . & f \\[6pt]
. & -b & . & . & f & -a & . & . & d \\
. & -c & . & . & h & . & -a & . & e \\
c & . & . & -h & . & . & . & -a & g
\end{matrix}
\label{eq!Tsyz}
\end{equation}
$M_1$ is of block form $(A\,B)$ with two $9\times8$ blocks, and one
checks that $LM_1=0$, and $M_1$ is isotropic for the standard quadratic
form $J=\left( \begin{smallmatrix} 0 & I \\ I & 0 \end{smallmatrix}
\right)$, so its kernel is $M_2=\left( \begin{smallmatrix} \tr B \\ \tr
A \end{smallmatrix} \right)$. The focus in \eqref{eq!Tsyz} is on $i$ as
an unprojection variable, multiplying $d,e,g,a$. One recognises its
Tom$_3$ matrix as the top $5\times5$ block, and the Koszul syzygy matrix
of $d,e,g,a$ as $\hbox{Submatrix}([6,7,8,14,15,16],[6,7,8,9])$; compare \cite{KM}.

For some of the spinors (see Section~\ref{s!3}), consider the
$8\times9$ submatrixes formed by 4 out of the first 5 rows of
\eqref{eq!Tsyz}, and the complementary 4 rows from the last 8. One
calculates their maximal minors with a mild effort:
\begin{equation}
\renewcommand{\arraycolsep}{0.3em}
\renewcommand{\arraystretch}{1.2}
\begin{array}{rcl}
\bigwedge^8M_{1,2,3,4,13,14,15,16} & = & a^2(af-be+cd)^2\cdot L, \\{}
\bigwedge^8M_{1,2,3,5,12,14,15,16} & = & a^2(ah-bg+ce)^2\cdot L, \\{}
\bigwedge^8M_{1,2,4,5,11,14,15,16} & = & a^2(-\la a^2-dg+e^2)^2\cdot L, \\{}
\bigwedge^8M_{1,3,4,5,10,14,15,16} & = & a^2(-\la ab-dh+ef)^2\cdot L, \\{}
\bigwedge^8M_{2,3,4,5,9,14,15,16} & = & a^2(-\la ac-eh+fg)^2\cdot L.
\end{array}
\label{eq!TsJ}
\end{equation}
The factor $a$ comes from the $3\times3$ diagonal block at the bottom
right, and the varying factors are the $4\times4$ Pfaffians of the first
$5\times5$ block. Compare \ref{s!Ks} for a sample Koszul syzygy.
\end{exa}

\begin{exc} \rm Apply column and isotropic row operations to put the
variable $f$ down a main diagonal of $B$; check that this puts the
complementary $A$ in the form of a skew $8\times8$ matrix and a row of
zeros. Hint: order the rows as $15,16,12,11,6,2,1,5,\
7,8,4,3,14,10,9,13$ and the columns as $1,2,-3,4,5,-7,8,9,6$. (See the
website for the easy code.) Do the same for either variable $e,h$, and
the same for any of $a,b,c$ (involving the multiplier $\la$).

Thus the isotropy condition $^t\!MJM$ can be thought of as {\em many}
skew symmetries.
\end{exc}

These examples provide useful sanity checks, with everything given by
transparent calculations; it is reassuring to be able to verify the
symmetry of the bilinear form on $P_2$ asserted in Proposition~1, the
shape of $A\,\tr B$ in \eqref{eq!AtB}, which parity of $J$ gives nonzero
spinors $\si_J$, and other minor issues of this nature.

I have written out the matrixes, spinors, Koszul syzygies etc.\ in a
small number of more complicated explicit examples (see the website). It
should be possible to treat fairly general Tom and Jerry constructions
in the same style, although so far I do not know how to use this to
predict anything useful. The motivation for this paper came in large
part from continuing attempts to understand Horikawa surfaces and Duncan
Dicks' 1988 thesis \cite{Di}, \cite{R1}.

\subsection{Definition of the Spin-Hom variety $\SpH_k$}\label{s!SpH}
Define the {\em Spin-Hom variety} $\SpH_k\subset\Mat(k+1,2k)$ as the
closure under $G=\GL(k+1)\times\Orth(2k)$ of the orbit of
$M^0=\left(\begin{smallmatrix} I_k & 0 \\ 0 & 0
\end{smallmatrix}\right)$, the second matrix in \eqref{eq!genf}. It
consists of isotropic homomorphisms $V_1\ot V_2$, in other words
matrixes $M_1$ whose $k+1$ rows are isotropic and mutually orthogonal
vectors in $V_2$ w.r.t.\ the quadratic form $\fie$, and span a subspace
that is in the given component of maximal isotropic subspaces if it is
$k$-dimensional.

In more detail, write $\SpH^0_k=G\cdot M^0\subset\Mat(k+1,2k)$ for the
orbit, $\SpH_k$ for its closure, and
$\SpH^{\dgn}_k=\SpH_k\setminus\SpH^0_k$ for the {\em degeneracy locus},
consisting of matrixes of rank $<k$. Section~\ref{s!3} discusses several
further properties of $\SpH_k$ and its degeneracy locus $\SpH^{\dgn}_k$.

\subsection{The Main Theorem} \label{s!main}
{\em Assume that $S$ is a polynomial ring graded in positive degrees.
Let $I$ be a
homo\-geneous ideal defining a codimension~$4$ Gorenstein subscheme
$X=V(I)\subset\Spec S$. Then a choice of minimal generators of $I$ (made
up of $k+1$ elements, say) and of the first syzygies between these
defines a morphism $\al\colon\Spec S\to \SpH_k$ such that
$\al\1(\SpH^{\dgn})$ has the same support as $X$, and hence codimension
$4$ in $\Spec S$.

Conversely, let $\al\colon\Spec S\to\SpH_k \subset \Mat(k+1,2k)$ be a
morphism for which $\al\1(\SpH^{\dgn})$ has codimension $\ge4$ in $\Spec
S$. Assume that $\al$ is graded, that is, equivariant for a positively
graded action of\/ $\Gm$ on $\SpH_k\subset\Mat(k+1,2k)$. Let
$M_1=(A\,B)$ be the matrix image of $\al$ (the matrix entries of $M_1$
or the coordinates of $\al$ are elements of $S$). Then by construction
$M_1$ and $J\,\tr M_1$ define the two middle morphisms of a complex. I
assert that this extends to a complex
\begin{equation}
0 \ot P_0 \xleftarrow{L} P_1 \xleftarrow{M_1} P_2 \xleftarrow{J\,\tr M_1}
P_3 \xleftarrow{\tr L} P_4 \ot 0. \label{eq!cx}
\end{equation}
in which $P_0,P_4\iso S$, the complex is exact except at $P_0$, and the
image of $L=(L_1,\dots,L_{k+1})$ generates the ideal of a Gorenstein
codimension~$4$ subscheme $X\subset\Spec S$. }

\subsection{Proof} \label{s!mpf}
The first part follows from what I have already said. The converse
follows by a straightforward application of the exactness criterion of
\cite{BE1}.

The complex $P\lodot$ of \eqref{eq!cx} comes directly from $M_1$.
Namely, define $P_0$ as the reflexive hull of $\coker\{M_1\colon P_1\ot
P_2$\} (that is, double dual); it has rank~1 because $M_1$ has generic
rank $k$. A graded reflexive module of rank~1 over a graded regular ring
is free (this is the same as saying that a Weil divisor on a nonsingular
variety is Cartier), so $P_0\iso S$. Given $P_3\iso P_1^\vee$, the
generically surjective map $S\iso P_0\ot P_1$ is dual to an inclusion
$S\into P_3$ that maps to the kernel of $P_2\ot P_3$.

The key point is to prove exactness of the complex
\[
P_0 \xleftarrow{\fie_1} P_1 \xleftarrow{\fie_2} P_2 \xleftarrow{\fie_3}
P_3 \xleftarrow{\fie_4} P_4 \ot 0,
\]
where I write $\fie_1=(L_1,\dots,L_{k+1})$, $\fie_2=M_1$, etc.\ to agree
with \cite{BE1}. The modules and homomorphisms
$P_0,\fie_1,P_1,\fie_2,P_2,\fie_3,P_3,\fie_4,P_4$ of this complex have
respective ranks $1,1,k+1,k,2k,k,k+1,1,1$, which accords with an exact
sequence of vector spaces, as in (\ref{eq!Vi}--\ref{eq!genf}); this is
Part (1) of the criterion of \cite{BE1}, Theorem~1.

The second condition Part (2) requires the matrixes of $\fie_i$ to have
maximal nonzero minors generating an ideal $I(\fie_i)$ that contains a
regular sequence of length $i$. However, $P\lodot$ is exact outside the
degeneracy locus, that is, at any point $\xi\in\Spec S$ for which
$\al(\xi)\notin\SpH_k^\dgn$, and by assumption, the locus of such points
has codimension $\ge4$. Thus the maximal minors of each $\fie_i$
generate an ideal defining a subscheme of codimension $\ge4$. In a
Cohen--Macaulay ring, an ideal defining a subscheme of codimension $\ge
i$ has height $\ge i$. \quad Q.E.D.

\section{Properties of $\SpH_k$ and its spinors} \label{s!3}

This section introduces the {\em spinors} as sections of the spinor line
bundle $\sS$ on $\SpH_k$. The {\em nonspinors} vanish on $\SpH_k$ and
cut it out in $V_k$ set theoretically. The spinors vanish on the {\em
other component} $\SpH_k'$ and cut out set theoretically the degeneracy
locus $\SpH^{\dgn}_k$ in $\SpH_k$.

The easy bit is to say that a spinor is the square root of a determinant
on $V_k\subset\Mat(k+1,2k)$ that vanishes to even order on a divisor of
$\SpH_k$ because it is locally the square of a Pfaffian. The ratio of
two spinors is a rational function on $\SpH_k$.

The tricky point is that the spinors are sections of the spinor bundle
$\sS$ on $\SpH_k$ that is defined as a $\Pin(2k)$ equivariant bundle, so
not described by any particularly straightforward linear or multilinear
algebra. As everyone knows, the spinor bundle $\sS$ on $\OGr(k,2k)$ is
the ample generator of $\Pic(\OGr(k,2k))$, with the property that
$\sS^{\tensor2}$ is the restriction of the Pl\"ucker bundle $\Oh(1)$ on
$\Gr(k,2k)$. On the affine orthogonal Grassmannian in Pl\"ucker space
$a\Gr(k,2k)\subset\bigwedge^k\kk^{2k}$, it corresponds to a 2-torsion
Weil divisor class. I write out a transparent treatment of the first
example in \ref{s!baby}.

I need to argue that the spinors pulled back to my regular ambient
$\Spec S$ by the appropriate birational transform are elements of $S$
(that is, poly\-nomials), rather than just sections of a spinor line
bundle. The reason that I expect to be able to do this is because I have
done many calculations like the Tom unprojection of \ref{s!exa},
Example~\ref{exa!T}, and it always works. In the final analysis, I win
for the banal reason that the ambient space $\Spec S$ has no 2-torsion
Weil divisors in its class group (because $S$ is factorial), so that the
birational transform of the spinor bundle $\sS$ to $\Spec S=\aff^N$ is
trivial.

The Cramer-spinor coordinates of the syzygy matrix $M_1=(A\,B)$ have the
potential to clarify many points about Gorenstein codimension~4: the
generic rank of $M_1$ is $k$, but it drops to $k-3$ on $\Spec R$; its
$k\times k$ minors are in $I_R^3$. There also seems to be a possible
explanation of the difference seen in examples between $k$ even and odd
in terms of the well known differences between the Weyl groups $D_k$
(compare \ref{s!evo}).

\subsection{Symmetry}
View $\GL(k+1)$ as acting on the first syzygy matrix $M_1(R)$ by row
operations, and $\Orth(2k)$ as column operations preserving the
orthogonal structure $\fie$, or the matrix $\J$. The maximal torus
$\Gm^{k+1}$ and Weyl group $A_k=S_{k+1}$ of the first factor
$\GL(k+1)$ act in the obvious way by scaling and permuting the rows of
$M_1$.

I need some standard notions for the symmetry of $\Orth(2k)$ and its
spinors. For further details, see Fulton and Harris \cite{FH}, esp.\
Chapter~20 and \cite{CR}, Section~4. Write $V_2=\kk^{2k}$ for the $2k$
dimensional vector space with basis $e_1,\dots,e_k$ and dual basis
$f_1,\dots,f_k$, making the quadratic form $\fie=\J$. Write
$U=U^k=\Span{e_1,\dots,e_k}$, so that $V_2=U\oplus U^\vee$. The
orthogonal Grassmannian $\OGr(k,2k)$ is defined as the variety of
$k$-dimensional isotropic subspaces that intersect $U$ in even
codimension, that is, in a subspace of dimension $\equiv k$ modulo~2.

\subsubsection{The $D_k$ symmetry of $\OGr(k,2k)$ and $\SpH_k$}
I describe the $D_k$ Weyl group symmetry of the columns in this notation
(compare \cite{CR}, Section~4). The maximal torus $\Gm^k$ of $\Orth(2k)$
multiplies $e_i$ by $\la_i$ and $f_i$ by $\la_i\1$, and acts likewise on
the columns of $M_1=(A\,B)$. The Weyl group $D_k$ acts on the $e_i,f_i$
and on the columns of $M_1=(A\,B)$ by permutations, as follows: the
subgroup $S_k$ permutes the $e_i$ simultaneously with the $f_i$; and the
rest of $D_k$ swaps {\em evenly many} of the $e_i$ with their
corresponding $f_i$, thus taking $U=\Span{e_1,\dots,e_k}$ to another
coordinate $k$-plane in $\OGr(k,2k)$. Exercise: The younger reader may
enjoy checking that the $k-1$ permutations
$s_i=(i,i+1)=(e_ie_{i+1})(f_if_{i+1})$ together with
$s_k=(e_kf_{k+1})(e_{k+1}f_k)$ are involutions satisfying the standard
Coxeter relations of type $D_k$, especially $(s_{k-1}s_k)^2=1$ and
$(s_{k-2}s_k)^3=1$.

\subsubsection{Spinor and nonspinor subsets} \label{s!sb}
The spinor sets $J\cup J^c$ index the spinors $\si_J$ (introduced in
\ref{s!sc}). Let $\{e_i,f_i\}$ be the standard basis of $\kk^{2k}$ with
form $\fie=\J$. There are $2^k$ choices of maximal isotropic subspaces
of $\kk^{2k}$ based by a subset of this basis; each is based by a subset
$J$ of $\{e_1,\dots,e_k\}$ together with the complementary subset $J^c$
of $\{f_1,\dots,f_k\}$. The spinor subsets are those for which $\#J$ has
the same parity as $k$, or in other words, the complement $\#J^c$ is
even; the nonspinor subsets are those for which $\#J$ has the parity of
$k-1$. The spinor set indexes a basis $\si_J$ of the spinor space of
$\OGr(k,2k)$, and similarly, the nonspinor set indexes the nonspinors
$\si'_{J'}$ of his dark twin $\OGr'(k,2k)$.

The standard affine piece of $\OGr(k,2k)$ consists of $k$-dimensional
spaces based by $k$ vectors that one writes as a matrix $(I\,A)$ with
$A$ a skew $k\times k$ matrix. The spinor coordinates of $(I\,A)$ are
the $2i\times2i$ diagonal Pfaffians of $A$ for $0\le i\le [\frac k2]$.
They correspond in an obvious way to the spinor sets just defined and
{\em they are the spinors} apart from the quibble about taking an
overall square root and what bundle they belong to.

\subsubsection{Even versus odd} \label{s!evo}
The distinction between $k$ even or odd is crucial for anything to do
with $\Orth(2k)$, $D_k$, spinors, Clifford algebras, etc. The spinor and
nonspinor sets correspond to taking a subset $J$ of $\{e_1,\dots,e_k\}$
and the complementary set $J^c$ of $\{f_1,\dots,f_k\}$. The $2^k$
choices correspond to the vertices of a $k$-cube. When $k$ is even this
is a bipartite graph; the spinors and nonspinors form the two parts. By
contrast, for odd $k$, both spinors and nonspinors are indexed by the
vertices of the $k$-cube divided by the antipodal involution (\cite{CR},
Section~4 writes out the case $k=5$ in detail).

For simplicity, I assume that $k$ is even in most of what follows; the
common case in applications that I really care about is $k=8$. Then
$J=\emptyset$ and $J^c=\{1,\dots,k\}$ is a spinor set, and the affine
pieces represented by $(I\,X)$ and $(Y\,I)$ (with skew $X$ or $Y$) are
in the same component of $\OGr(k,2k)$. The odd case involves related
tricks, but with some notable differences of detail (compare \cite{CR},
Section~4).

\subsubsection{The other component $\OGr'$ and $\SpH_k'$}

I write $\OGr'(k,2k)$ for the {\em other component} of the maximal
isotropic Grassmannian, consisting of subspaces meeting $U$ in odd
codimension. Swapping {\em oddly many} of the $e_i$ and $f_i$
interchanges $\OGr$ and $\OGr'$. Likewise, $\SpH_k'$ is the closure of
the $G$-orbit of the matrix $M_0'$ obtained by interchanging one
corresponding pair of columns of $M_0$.

\begin{cla} \label{cla!Vk}
Write $V_k$ for the scheme defined by \eqref{eq!AtB} (that is, the
``variety of complexes''). It has two irreducible components
$V_k=\SpH_k\cup\SpH'_k$ containing matrixes of maximal rank $k$. The two
components are generically reduced and intersect in the degenerate locus
$\SpH^{\dgn}_k$. (But one expects $V_k$ to have embedded primes at its
smaller strata, as in the discussion around \eqref{eq!sq}.)
\end{cla}

This follows from the properties of spinor minors $\De_J$ discussed in
Exercise~\ref{exc!s}: the $\De_J$ are $k\times k$ minors defined as
polynomials on $V_k$, and vanish on $\SpH'_k$ but are nonzero on a dense
open subset of $\SpH_k$.

\subsection{A first introduction to $\OGr(k,2k)$ and its spinors}
\label{s!baby}

The lines on the quadric surface provide the simplest calculation, and
already have lots to teach us about $\OGr(2,4)$ and $\OGr(k,2k)$: the
conditions for the $2\times4$ matrix
\begin{equation}
N=\begin{pmatrix} a & b\; & x & y \\ c & d\; & z & t \end{pmatrix}
\end{equation}
to be isotropic for $\J$ are
\begin{equation}
ax + by=0, \quad az + bt + cx + dy=0, \quad cz + dt=0.
\label{eq!2x4}
\end{equation}
Three equations \eqref{eq!2x4} generate an ideal $I_W$ defining
a codimension~3 complete inter\-section $W\subset\aff^8$ that breaks up
into two components $\Si\sqcup\Si'$, corresponding to the two pencils of
lines on the quadric surface: the two affine pieces of $\OGr(2,4)$ that
consist of matrixes row equivalent to $(I\,A)$ or $(A\,I)$, with $A$ a
skew matrix, have one of the {\em spinor minors} $\De_1=ad-bc$ or
$\De_2=xt-yz$ nonzero, and
\begin{equation}
dx-bz=at-cy=0 \quad\hbox{and}\quad dy-bt=-(az-cx)
\end{equation}
on them. This follows because all the products of $\De_1,\De_2$ with the
{\em nonspinors minors} $dx-bz,at-cy$ are in $I_W$, as one checks
readily. Thus if $\De_1\ne0$ (say), I can multiply by the adjoint of the
first block to get
\begin{equation}
\begin{pmatrix} d & -b \\ -c &
a \end{pmatrix}
\begin{pmatrix} a & b & x & y \\ c & d & z & t \end{pmatrix} =
\begin{pmatrix} \De_1 & 0 & dx-bz & dy-bt \\ 0 & \De_1 & az-cx & at-cy
\end{pmatrix}
\label{eq!De}
\end{equation}
where the second block is skew. Note that
\begin{equation}
\De_1\cdot\bigl(\De_1\De_2-(az-cx)^2\bigr) \in I_W.
\label{eq!sq}
\end{equation}
If $\De_1\ne0$, the relations \eqref{eq!2x4} imply that we are in $\Si$.
The ideal of $\Si$ is obtained from \eqref{eq!2x4} allowing cancellation
of $\De_1$; in other words $I_\Si=[I_W:\De_1]$ is the colon ideal with
either of the spinor minors $\De_1$ or $\De_2$.

The second block in \eqref{eq!De} is only skew mod $I_W$ after
cancelling one of $a,b,\dots,t$; similarly $\De_1\De_2-(az-cx)^2\notin
I_W$, so that \eqref{eq!sq} involves cancelling $\De_1$. Thus a
geometric description of $\Si,\Si'\subset\Mat(k,2k)$ should usually lead
to ideals with embedded primes at their intersection or its smaller
strata.

Now by relation \eqref{eq!sq}, the Pl\"ucker embedding takes $\OGr(2,4)$
to the conic $XZ=Y^2$, with $X=\De_1=ad-bc$, $Y=az-cx$, $Z=\De_2=xt-yz$.
This is $(\PP^1,\Oh(2))$ parametrised by $u^2,uv,v^2$ where $u,v$ base
$H^0(\PP^1,\Oh(1))$. Thus $X=u^2$, $Y=uv$ and $Z=v^2$ on $\OGr(2,4)$;
the spinors are $u$ and $v$. The ratio $u:v$ equals $X:Y=Y:Z$. Each of
$\De_1$ and $\De_2$ vanishes on a double divisor, but the quantities
$u=\sqrt{\De_1}$, $v=\sqrt{\De_2}$ are not themselves polynomial.

The conclusion is that the minors $\De_1$ and $\De_2$ are {\em spinor
squares}, that is, squares of sections $u,v$ of a line bundle $\sS$, the
spinor bundle on $\OGr(2,4)$. If we view $\OGr(2,4)$ as a subvariety of
$\Gr(2,4)$, only $\sS^{\tensor2}$ extends to the Pl\"ucker line bundle
$\Oh(1)$. Embedding $\OGr(2,4)$ in the Pl\"ucker space
$\PP(\bigwedge^2\C^4)$ and taking the affine cone gives the affine
spinor variety $a\OGr(2,4)$ as the cone over the conic, and $\sS$ with
its sections $u,v$ as the ruling.

In fact $a\OGr(2,4)$ and his dark twin $a\OGr'$ are two ordinary quadric
cones in linearly disjoint vector subspaces of the Pl\"ucker space
$\bigwedge^2\C^4$, and the spinor bundle on the union has a divisor
class that is a 2-torsion Weil divisor on each component. This picture
is of course the orbifold quotient of $\pm1$ acting on two planes
$\aff^2$ meeting transversally in $\aff^4$.

\subsubsection{Exercise} \label{exc!s}
Generalise the above baby calculation to the subvariety
$W_k\subset\Mat(k,2k)$ of matrixes $(A\,X)$ whose $k$ rows span an
isotropic space for $\J$, or in equations, the $k\times k$ product
$A\,\tr X$ is skew. Assume $k$ is even.

\begin{enumerate}
\renewcommand{\labelenumi}{(\arabic{enumi})}

\item $W_k\subset\Mat(k,2k)$ is a complete intersection subvariety of
codimension $\binom{k+1}2$. [Hint: Just a dimension count.]

\item $W_k$ breaks up into two irreducible components $\Si\cup\Si'$,
where $\Si$ contains the space spanned by $(I\,X)$ with $X$ skew, or
more generally, by the span of the columns $J\cup J^c$ for $J$ a spinor
set; its nondegenerate points form a principal $\GL(k)$ bundle over the
two components $\OGr\sqcup\OGr'$ of the maximal isotropic Grassmannian.

\item For $J$ a spinor set, the $k\times k$ spinor minor $\De_J$ of
$(A\,X)$ (the determinant of the submatrix formed by the columns $J\cup
J^c$) is a polynomial on $\Mat(k\times2k)$ that vanishes on $\Si'$, and
vanishes along a double divisor of $\Si$, that is, twice a prime Weil
divisor $D_J$.

\item The Weil divisors $D_{J_1}$ and $D_{J_2}$ corresponding to two
spinor sets $J_1$ and $J_2$ are linearly equivalent. [Hint: First
suppose that $J_1$ is obtained from $J$ by exactly two transpositions,
say $(e_1f_2)(e_2f_1)$, and argue as in \eqref{eq!sq} to prove that
$\si_J\si_{J_1}$ restricted to $\Si$ is the square of either minor
obtained by just one of the transpositions.]

\end{enumerate}

\subsubsection{Spinors on $\OGr(k,2k)$} \label{ss!stuff}
The orthogonal Grassmann variety $\OGr(k,2k)$ has a {\em spinor}
embedding into $\PP(\kk^{2^{k-1}})$, of which the usual Pl\"ucker
embedding
\[
\OGr(k,2k)\subset\Gr(k,2k) \into \PP\Bigl(\bigwedge^k\kk^{2k}\Bigr)
\]
is the Veronese square. The space of spinors $\kk^{2^{k-1}}$ is a
representation of the spin double cover $\Pin(2k)\to\Orth(2k)$.

A point $W\in\OGr(k,2k)$ is a $k$-dimensional subspace $W^k\subset
\kk^{2k}$ isotropic for $\J$ and intersecting $U=\Span{e_1,\dots,e_n}$
in even codimension. I can write a basis as the rows of a $k\times 2k$
matrix $N_W$. If I view $W$ as a point of $\Gr(k,2k)$, its Pl\"ucker
coordinates are all the $k\times k$ minors of $N_W$. There are
$\binom{2k}k$ of these (that is, 12870 if $k=8$), a fraction of which
vanish $\OGr(k,2k)$, as the determinant of a skew matrix of odd size.

The finer embedding of $\OGr(k,2k)$ is by spinors. The spinors $\si_J$
are sections of the spinor line bundle $\sS$, $2^{k-1}$ of them (which
is $128$ if $k=8$, about $1/100$ of the number of Pl\"ucker minors).
Each comes by taking a $k\times k$ submatrix formed by a spinor subset
of columns of $N_W$ (in other words, restricting to an isotropic
coordinate subspace of $\kk^{2k}$ in the specified component
$\OGr(k,2k)$), taking its $2\ka\times2\ka$ minor (where $\ka=\left[\frac
k2\right]$) and factoring it as the perfect square of a section of
$\sS$. The only general reason for a $2\ka\times2\ka$ minor to be a
perfect square is that the submatrix is skew in some basis; in fact, as
in \eqref{eq!De}, after taking one fixed square root of a determinant,
and making a change of basis, the maximal isotropic space can be written
as $(I\,X)$ with $X$ skew, and the spinors are all the Pfaffians of $X$.

\subsection{Cramer-spinor coordinates on $\SpH_k$} \label{s!sc}
\subsubsection{Geometric interpretation}
A point of the open orbit $\SpH_k^0\subset\SpH_k$ is a matrix $M$ of
rank $k$; it defines an isomorphism from a $k$-dimensional subspace of
$V_1$ (the column span of $M$) to its row span, a maximal isotropic
subspace of $V_2$ in the specified component $\OGr(k,2k)$. Therefore the
nondegenerate orbit $\SpH_k^0\subset\SpH_k$ has a morphism to
$\PP(V_1^\vee)\times\OGr(k,2k)$ that makes it a principal $\GL(k)$
bundle. The product $\PP(V_1^\vee)\times\OGr(k,2k)$ is a projective
homogeneous space under $G=\GL(k+1)\times\Pin(2k)$

It embeds naturally in the projectivisation of
$\kk^{k+1}\tensor\kk^{2^{k-1}}\!$, with the second factor the space of
spinors. This is the representation of $G$ with highest weight vector
$v=(0,\dots,0,1)\tensor(1,0,\dots,0)$. The composite
\begin{equation}
\SpH_k^0 \to \PP(V_1^\vee)\times\OGr(k,2k) \into
\PP(\kk^{k+1}\tensor\kk^{2^{k-1}})
\label{eq!Csp}
\end{equation}
takes the typical matrix $M_0$ (or equivalently, the complex
\eqref{eq!genf}) to $v$.

The Cramer-spinor coordinates of $\al\in\SpH_k(S)$ are the bihomogeneous
coordinates under the composite map \eqref{eq!Csp}.

\subsubsection{Spinors as polynomials} \label{s!spin}

The spinors $\si_J$ occur naturally as sections of the spinor line
bundle $\sS$ on $\OGr(k,2k)$, and so have well defined pullbacks to
$\SpH_k^0$ or to any scheme $T$ with a morphism $\al\colon
T\to\SpH_k^0$. For $\si_J$ to be well defined in $H^0(\Oh_T)$, the
pullback of the spinor line bundle to $T$ must be trivial.

\begin{lem}
Let $\al\in\Mor(\Spec S,\SpH_k)=\SpH_k(S)$ be a classifying
map as in Theorem~\ref{s!main} and write $M_1\in \Mat(S,k+1,2k)$ for its
matrix (with entries in $S$). Then for a spinor set $J\cup J^c$ (as in
\ref{s!sb}), the $(k+1)\times k$ submatrix $N_J$ of $M_1$ with columns
$J\cup J^c$ has
\begin{equation}
\bigwedge\nolimits^k N_J= L\cdot\si_J^2,
\label{eq!wdg}
\end{equation}
where $L=(L_1,\dots,L_{k+1})$ generates the cokernel of $M_1$, and
$\si_J\in S$.
\end{lem}

\subsection{Proof} A classifying map $\al\in\SpH_k(S)$ as in
Theorem~\ref{s!main} restricts to a morphism $\al$ from the
nondegenerate locus $\Spec S\setminus V(I_R)$ to $\SpH_k^0$; on the
complement of $V(I_R)$, the matrix $M_1$ has rank $k$, and its $k$th
wedge defines the composite morphism to the product
$\PP^k\times\Gr(k,2k)$ in its Segre embedding:
\begin{multline}
\Spec S\setminus V(I_R) \to \SpH^0_k\to\PP^k\times\OGr(k,2k) \\
\into \PP^k\times\Gr(k,2k) \subset\PP\left(\kk^{k+1}\tensor
\bigwedge\nolimits^kV^{2k}\right).
\end{multline}
The entries of $\bigwedge^k N_J$ are $k+1$ coordinates of this morphism,
and are of the form $L_i\cdot\si_J^2$ already on the level of
$\PP^k\times\OGr(k,2k)$.

Note that $\Spec S\setminus V(I_R)$ is the complement in $\Spec
S=\aff^N$ of a subset of codimension $\ge4$ so has trivial $\Pic$. Each
maximal minor of $N_J$ splits as $L_i$ times a polynomial that vanishes
on a divisor that is a double (because it is the pullback of the square
of a spinor); therefore the polynomial is a perfect square in $S$.
\QED\par\medskip

The following statement is the remaining basic issue that I am currently
unable to settle in general.

\begin{conj}
Under the assumptions of Lemma~\ref{s!spin}, $\si_J\in I_R$.
\end{conj}

This is clear when $R$ is reduced, that is, $I_R$ is a radical ideal.
Indeed if $\si_J$ is a unit at some generic point $\xi\in V(I_R)=\Spec
R$, then \eqref{eq!wdg} implies that $I_R$ is generated at $\xi$ by the
$k\times k$ minors of the $(k+1)\times k$ matrix $N_J$; these equations
define a codimension~2 subscheme of $\Spec S$, which is a contradiction.
This case is sufficient for applications to construction of ordinary
varieties, but not of course to Artinian subschemes of $\aff^4$.

The conjecture also holds under the assumption that $I_R$ is generically
a codimension~4 complete intersection. Indeed, the resolution of $I_R$
near any generic point $\xi\in V(I_R)$ is then the $4\times6$ Koszul
resolution of the complete intersection direct sum some nonminimal stuff
that just add invertible square matrix blocks. Then both the $L_i$ and
the $\si_J$ are locally given by Example~\ref{exa!Kos4}.

At present, the thing that seems to make the conjecture hard is that the
definition of the $\si_J$ and the methods currently available for
getting formulas for them consists of working on the nondegenerate locus
of $\SpH_k$: choose a block diagonal form and take the Pfaffian of a
skew complement, \dots This is just not applicable at points $\si\in
V(I_R)$.

The conjecture could possibly be treated by a more direct understanding
of the spin morphism $\Spec S\to\kk^{2k}$ defined by spinors and
nonspinors, not passing via the square root of the Pl\"ucker morphism as
I do implicitly in Lemma~1 by taking $\bigwedge^k$.

\section{Final remarks, open problems} \label{s!4}

\subsection{Birational structure and dimension of $\SpH_k$} A general
$M=(A\,B)\in\SpH_k$ has $k+1$ rows that span a maximal isotropic space
$U\in\OGr(k,2k)$ and $2k$ columns that span a $k$-dimensional vector
subspace of $\kk^{k+1}$, that I can view as a point of $\PP^k$; thus
$\SpH_k^0$ is a principal $\GL(k)$ bundle over $\PP^k\times\OGr(k,2k)$.
In particular, $\dim\SpH_k=k^2+k+\binom k2=\frac{3k^2+k}2$.

The tangent space to $\SpH_k$ at the general point
$M_0=\left(\begin{smallmatrix} I_k & 0 \\ 0 & 0
\end{smallmatrix}\right)$ is calculated by writing an infinitely near
matrix as $M_0+\left(\begin{smallmatrix} A'_k & B'_k \\ a_{k+1} &
b_{k+1} \end{smallmatrix}\right)$; here the blocks $A'_k$ and $B'_k$ are
$k\times k$ matrixes, and $a_{k+1}$ and $b_{k+1}$ are $1\times k$ rows.
Then the tangent space to $V_k$ defined by $A\,\tr B=0$ is the affine
subspace obtained by setting $B'_k$ to be skew and $b_{k+1}=0$.
Therefore $\SpH_k$ has codimension $\binom{k+1}2+k$ and dimension
$2k(k+1)-\binom{k+1}2-k=\frac{3k^2+k}2$.

It is interesting to observe that equations \eqref{eq!AtB} express
$\SpH_k\cup\SpH'_k$ as an almost complete intersection. Namely,
\eqref{eq!AtB} is a set of $\binom{k+1}2$ equations in $\aff^{2k(k+1)}$
vanishing on a variety of dimension $\frac{3k^2+k}2$, that is, of
codimension $\binom{k+1}2-1$.

\subsection{Intermediate rank} \label{s!ir}

The Spin-Hom variety $\SpH_k$ certainly contains degenerate matrixes
$M_1$ of rank $k-1$ or $k-2$, but any morphism $\Spec S\to\SpH_k$ that
hits one of these must hit the degeneracy locus in codimension $\le3$,
so does not correspond to anything I need here. The following claim must
be true, but I am not sure where it fits in the logical development.

\begin{cla} Every point $P\in\SpH_k$ corresponds to a matrix
$M_1=(A\,B)$ of rank $\le k$. If a morphism $\al\colon\Spec S\to\SpH_k$
takes $\xi$ to a matrix $M_1$ of rank $k+1-i$ for $i = 1,2,3,4$ then
$\al\1(\SpH_k^\dgn)$ has codimension $\le i$ in a neighbourhood of
$\xi$. In other words, a morphism $\al$ that is regular in the sense of
my requirement never hits matrixes $M_1$ of rank intermediate between
$k$ and $k-3$; and if $\al$ is regular then $\al\1(\SpH_k^\dgn)$ has
codimension exactly~$4$. \end{cla}

\subsection{The degeneracy locus as universal subscheme} \label{s!sbt}

The proof in \ref{s!mpf} doesn't work for $\SpH_k$ itself in a
neighbourhood of a point of $\SpH^{\dgn}_k$, because taking the
reflexive hull, and asserting that $P_0$ is locally free works only over
a regular scheme. Moreover, it is not just the proof that goes wrong. I
don't know what happens over the strata of $\SpH^{\dgn}_k$ where $M_1$
drops rank by only~1 or~2.

We discuss the speculative hope that $\SpH^{\dgn}_k\subset\SpH_k$ has a
description as a kind of universal codimension~4 subscheme, with the
inclusions enjoying some kind of Gorenstein adjunction properties. But
if this is to be possible at all, we must first discard uninteresting
components of $\SpH^{\dgn}_k$ corresponding to matrixes of intermediate
rank $k-1$ or $k-2$.

It is possible that there is some universal blowup of some big open in
$\SpH_k$ that supports a Gorenstein codimension~4 subscheme and would be
a universal space in a more conventional sense. Or, as the referee
suggests, there might be a more basic sense in which appropriate
codimension~4 components $\Ga$ of the degeneracy locus are {\em
universal Gorenstein embeddings}, meaning that the adjunction
calculation $\om_\Ga=\Ext^4_{\Oh_{\SpH}}(\Oh_\Ga,\om_{\SpH})$ for the
dualising sheaf is locally free and commutes with regular pullbacks.

\subsection{Koszul syzygies} \label{s!Ks}

Expressing the generators of $I$ as a function of the entries of the
syzygy matrix is essentially given by the map $\bigwedge^2 P_1\to P_2$
that writes the Koszul syzygies as linear combinations of the minimal
syzygies.

The $L_i$ are certainly linear combinations of the entries of $M_1$.
More precisely, since the $2k$ columns of $M_1$ provide a minimal basis
for the syzygies, they cover in particular the Koszul syzygies $L_i\cdot
L_j-L_j\cdot L_i\equiv0$. This means that for every $i\ne j$ there is
column vector $v_{ij}$ with entries in $S$ such that
$M_1v_{ij}=(\dots,L_j,\dots, L_i,\dots)$ is the column vector with $L_j$
in the $i$th place and $L_i$ in the $j$th and 0 elsewhere. For example,
referring to Example~\ref{exa!T}, you might enjoy the little exercise in
linear algebra of finding the vector
\begin{multline*}
v=(-\la c,\la b,0,0,0,d,e,g,0,0,0,0,0,0,0,0) \quad\hbox{for which} \\
v\, \tr M_1=(-\la ab-dh+ef,-\la ac-eh+fg,0,0,0,0,0,0,0),
\end{multline*}
where $\tr M_1$ is the matrix of \eqref{eq!Tsyz}, and similarly for 35
other values of $i,j$.

\subsection{More general ambient ring $S$} \label{ss!J}
I restrict to the case of ideals in a graded polynomial ring over a
field of characteristic $\ne2$ in the belief that progress in this case
will surely be followed by the more general case of a regular local
ring. Then $P_2$ is still a free module, with a perfect symmetric
bilinear form $S^2(P_2)\to P_4$, with respect to which $P_1\ot P_2$ is
the dual of $P_2\ot P_3$. This can be put in the form $\J$ over the
residue field $\kk_0=S/m_S$ of $S$ if we assume that $k(S)$ is
algebraically closed and contains $\frac12$; we can do the same over $S$
itself if we assume that $S$ is complete (to use Hensel's Lemma). At
some point if we feel the need for general regular rings, we can
probably live with a perfect quadratic form $\fie$ and the dualities it
provides, without the need for the normal form $\J$.

\subsection{More general rings and modules}
Beyond the narrow question of Gorenstein codimension~4, one could ask
for the structure of any free resolution of an $S$-module $M$ or
$S$-algebra $R$. As in \ref{s!gf}, one can say exactly what the general
fibre is, and think of the complex $P\lodot$ as a fibre bundle over
$S\setminus\Supp M$ with some product of linear groups as structure
group. If we are doing $R$-algebras, the complex $P\lodot$ also has a
symmetric bilinear structure, that reduces the structure group. My point
is that if we eventually succeed in making some progress with Gorenstein
codimension~4 rings, we might hope to also get some ideas about
Cohen--Macaulay codimension~3 and Gorenstein codimension~5.

For example, in vague terms, there is a fairly clear strategy how to
find a key variety for the resolution complexes of Gorenstein
codimension~5 ideals, by analogy with my Main Theorem~\ref{s!main}. In
this case, the resolution has the shape
\begin{equation}
0 \ot P_0 \ot P_1 \ot P_2 \ot P_3 \ot P_4 \ot P_5 \ot 0,
\label{eq!P5}
\end{equation}
with $P_0=S$, $P_1=(a+1)S$, $P_2=(a+b)S$ and $P_3,\dots,P_5$ their
duals. The complex is determined by two syzygy matrixes
$M_1\in\Mat(a+1,a+b)$ of generic rank $a$ defining $P_1 \ot P_2$ and a
symmetric $(a+b)\times(a+b)$ matrix $M_2$ of generic rank $b$ defining
$P_2 \ot P_3=P_2^\vee$, constrained by the complex condition $M_1M_2=0$.
The ``general fibre'' is given by the pair
$M_1=\left(\begin{smallmatrix} I_a & 0 \\ 0 & 0
\end{smallmatrix}\right)$, $M_2=\left(\begin{smallmatrix} 0 & 0 \\ 0 &
I_b \end{smallmatrix}\right)$, the appropriate key variety is its closed
orbit under $\GL(a+1)\times\GL(a+b)$. The maximal nonzero minors of
$M_1$ and $M_2$ define a map to a highest weight orbit in
\[
\Hom \Bigl(\bigwedge^a P_2, \bigwedge^a P_1\Bigr)
\times \Sym^2\Bigl(\bigwedge^b P_2\Bigr).
\]

\subsection{Difficulties with applications} \label{s!app}
I expand what the introduction said about the theory currently not being
applicable. We now possess hundreds of constructions of codimension~4
Gorenstein varieties, for example, the Fano 3-folds of \cite{TJ}, but
their treatment (for example, as Kustin--Miller unprojections) has
almost nothing to do with the structure theory developed here. My Main
Theorem~\ref{s!main} does not as it stands construct anything, because
it does not say how to produce morphisms $\al\colon\Spec S\to\SpH_k$, or
predict their properties. The point that must be understood is not the
key variety $\SpH_k$ itself, but rather the space of morphisms
$\Mor(\Spec S,\SpH_k)$, which may be intractable or infinitely
complicated (in the sense of Vakil's Murphy's law \cite{Va}); there are
a number of basic questions here that I do not yet understand.

Even given $\al$, we do not really know how to write out the equations
$(L_1,\dots,L_{k+1})$, other than by the implicit procedure of taking
hcfs of $k\times k$ minors. One hopes for a simple formula for the
defining relations $L_i$ as a function of the first syzygy matrix
$M_1=(A\,B)$. Instead, one gets the vector $(L_1,\dots,L_{k+1})$ by
taking out the highest common factor from $\bigwedge^k M_I$ for any
spinor subset $I$, asserting that it is a perfect square $\si_J^2$. The
disadvantage is that as it stands this is only implicitly a formula for
the $L_i$.

\subsection{Obstructed constructions}

One reason that $\Mor(S,\SpH_k)$ is complicated is that the target is
big and singular and needs many equations. However, there are also
contexts in which $S$-valued points of much simpler varieties already
give families of Gorenstein codimension~4 ideals that are obstructed in
interesting ways.

Given a $2\times4$ matrix $A=(\begin{smallmatrix} a_1&a_2&a_3&a_4 \\
b_1&b_2&b_3&b_4 \end{smallmatrix})$ with entries in a regular ring $S$,
the 6 equations $\bigwedge^2A=0$ define a Cohen--Macaulay codimension~3
subvariety $V\subset\Spec S$. An elephant $X\in|{-}K_V|$ is then a
Gorenstein subvariety of codimension 4 with a $9\times16$ resolution. If
we are in the ``generic'' case with 8 independent indeterminate entries,
$V$ is the affine cone over $\Segre(\PP^1\times\PP^3)$, and $X$ is a
cone over a divisor of bidegree $(k,k+2)$ in $\Segre(\PP^1\times\PP^3)$.

Although $X\subset V$ is a divisor, if we are obliged to treat it by
equations in the ambient space $\Spec S$, it needs 3 equations in
``rolling factors format''. The general case of this is contained in
Dicks' thesis \cite{Di}, \cite{R1}: choose two vectors $m_1,m_2,m_3,m_4$
and $n_1,n_2,n_3,n_4$, and assume that the identity
\begin{equation}
\sum a_in_i \equiv \sum b_im_i
\label{eq!aini}
\end{equation}
holds as an equality in the ambient ring $S$. Then the 3 equations
\begin{equation}
\sum a_im_i = \sum b_im_i \equiv \sum a_in_i = \sum b_in_i = 0
\label{eq!aimi}
\end{equation}
define a hypersurface $X\subset V$ that is an elephant $X\in|{-}K_V|$
and thus a Gorenstein subvariety with $9\times16$ resolution.

The problem in setting up the data defining $X$ is then to find
solutions in $S$ of \eqref{eq!aini}. In other words, these are
$S$-valued points of the affine quadric cone $Q_{16}$, or morphisms
$\Spec S\to Q_{16}$. How to map a regular ambient space to the quadratic
cone $Q_{16}$ is a small foretaste of the more general problem of the
classifying map $\Spec S\to\SpH_k$. This case is discussed further in
\cite{Ki}, Example~10.8, which in particular writes out explicitly the
relation between \eqref{eq!aimi} and the classifying map $\Spec
S\to\SpH_k$ of Theorem~\ref{s!main}.

There are many quite different families of solutions to this problem,
depending on what assumptions we make about the graded ring $S$, and how
general we take the matrix $A$ to be; different solutions have a number
of important applications to construction and moduli of algebraic
varieties, including my treatment of the Horikawa quintic $n$-folds.

Another illustration of the phenomenon arises in a recent preprint of
Catanese, Liu and Pignatelli \cite{CLP}. Take the $5\times5$ skew matrix
\begin{equation}
M=\begin{pmatrix}
v & u & z_2 & D \\
& z_1 & y & m_{25} \\
&& l & m_{35} \\
&&& m_{45}
\end{pmatrix}
\label{eq!FRW}
\end{equation}
with entries in a regular ring $S_0$, and suppose that $v,u,z_2,D$ forms
a regular sequence in $S$. Assume that the identity
\begin{equation}
z_1m_{45} - ym_{35} + lm_{25} \equiv av + bu + cz_2 + dD
\label{eq!Pf2345}
\end{equation}
holds as an equality in $S_0$. The identity \eqref{eq!Pf2345} puts the
Pfaffian $\Pf_{23.45}$ in the ideal $(v,u,z_2,D)$; the other 4 Pfaffians
are in the same ideal for the trivial reason that every term involves
one entry from the top row of $M$.

This is a new way of setting up the data for a Kustin--Miller
unprojection: write $Y\subset\Spec S_0$ for the codimension~3 Gorenstein
subscheme defined by the Pfaffians of $M$. It contains the codimension~4
complete intersection $V(v,u,z_2,D)$ as a codimension~1 subscheme, and
unprojecting $V$ in $Y$ adjoins an unprojection variable $x_2$ having 4
linear equations $x_2\cdot(v,u,z_2,D)=\cdots$, giving a codimension~4
Gorenstein ring with $9\times16$ resolution.

The problem of how to fix \eqref{eq!Pf2345} as an identity in $S_0$ is
again a question of the $S_0$-valued points of a quadric cone, this time
a quadric $Q_{14}$ of rank~14. \cite{CLP}, Proposition~5.13 find two
different families of solutions, and exploit this to give a local
description of the moduli of their surfaces.

At first sight this looks a bit like a Jerry$_{15}$ unprojection. In
fact one of the families of \cite{CLP} (the one with $c_0=B_x=0$) can
easily be massaged to a conventional Jerry$_{15}$ having a double Jerry
structure (compare \cite{TJ}, 9.2), but this does not seem possible for
the more interesting family in \cite{CLP} with $D_x=(l/c_0)B_x$.

\medskip\noindent{\bf Question} Do these theoretical calculations
contain the results of \cite{Di}, \cite{CLP} and the like?

\medskip\noindent{\bf Answer} Absolutely not. They may provide a
framework that can produce examples, or simplify and organise the
construction of examples. To get complete moduli spaces, it is almost
always essential to use other methods, notably infinitesimal deformation
calculations or geometric constructions.

\medskip\noindent{\bf Question} The fact that $S$ can have various gradings seems to add
to the complexity of the space $\Mor(S,\SpH_k)$, doesn't it?

\medskip\noindent{\bf Answer} That may not be the right interpretation
-- we could perhaps think that $\Mor(S,\SpH_k)$ (or even the same just
for $\Mor(S,Q_{2k})$ into a quadric of rank $2k\ge4$) is infinite
dimensional and infinitely complicated, so subject to Murphy's law
\cite{Va}, but that when we cut it down to graded in given degrees, it
becomes finitely determined, breaking up into a number of finite
dimensional families that may be a bit singular, but can be studied with
success in favourable cases.

\subsection{Problem session}

\subsubsection{Computing project} It is a little project in computer
algebra to write an algorithm to put the projective resolution
\eqref{eq!P4} in symmetric form. This might just be a straightforward
implementation of the Buchsbaum--Eisenbud symmetrised complex
$S^2P\lodot$ outlined in Section~\ref{s!1}. Any old computer algebra
package can do syzygies, but as far as I know, none knows about the
symmetry in the Gorenstein case.

We now have very many substantial working constructions of codimension~4
Gorenstein varieties. We know in principle that the matrix of first
syzygies can be written out in the $(A\,B)$ form of \eqref{eq!Tsyz},
but as things stand, it takes a few hours or days of pleasurable
puzzling to do any particular case.

\subsubsection{Linear subvarieties} \label{s!lin}
What are the linear subvarieties of $\SpH_k$? The linear question may be
tractable, and may provide a partial answer to the quest for an explicit
structure result.

The Spin-Hom variety $\SpH_k$ is defined near a general point by
quadratic equations, so its linear subspaces can be studied by the
tangent-cone construction by analogy with the linear subspaces of
quadrics, Segre products or Grassmannians: the tangent plane $T_P$ at
$P\in V$ intersects $V$ in a cone, so that linear subspaces of $V$
through $P$ correspond to linear subspaces in the base of the cone. Now
choose a point of the projected variety and continue.

Presumably at each stage there are a finite number of strata of the
variety in which to choose our point $P$, giving a finite number of
types of $\Pi$ up to symmetry. I believe that the two famous cases of
the Segre models of $\PP^2\times\PP^2$ and $\PP^1\times\PP^1\times\PP^1$
are maximal linear space of $\SpH_8$.

It is possible that this method can be used to understand more general
morphisms $\Spec S\to\SpH_k$ from the regular space $\Spec S$. In this
context, it is very suggestive that Tom and Jerry \cite{TJ} are given in
terms of linear subspaces of $\Gr(2,5)$. In this case, the intersection
with a tangent space is a cone over $\PP^1\times\PP^2$, so it is clear
how to construct all linear subspaces of $\Gr(2,5)$, and equally clear
that there are two different families, and how they differ.

\subsubsection{Breaking the $A_k$ and $D_k$ symmetry} \label{s!break}

Experience shows that the bulk constructions of Gorenstein codimension~4
ideals do not have the symmetry of the Buchsbaum--Eisenbud Pfaffians in
codimension~3. The equations and syzygies invariably divide up into
subsets that one is supposed to treat inhomogeneously. For example, in
the $9\times16$ unprojection cases, the defining equations split into
two sets, the 5 Pfaffian equations of the variety in codimension~3 not
involving the unprojection variable $s$, and the 4 unprojection
equations that are linear in $s$.

The columns of the syzygy matrix $(A\,B)$ are governed by the algebraic
group $\Spin(2k)$ of type $D_k$, whereas its rows are governed by
$\GL(k+1)$ of type $A_k$. The common bulk constructions of Gorenstein
codimension~4 ideals seem to to accommodate the $A_k$ symmetry of the
rows of $M_1$ and the $D_k$ symmetry of its columns by somehow breaking
both to make them compatible. This arises if you try to write the 128
spinor coordinates $\si_J$ as linear combinations of the 9 relations
$(L_1,\dots,L_{k+1})$, so relating something to do with the columns of
$M_1$ to its rows. This symmetry breaking and its effect is fairly
transparent in \ref{s!exa}, Example~\ref{exa!PfHy}, \eqref{eq!Pfh}.

Example~\ref{exa!T} is more typical. (This case comes with three
different Tom projections, so may be more amenable.)  Of the 128 spinors
$\si_J$, it turns out that 14 are zero, 62 are of the form a monomial
times one of the relations $L_i$ (as in \eqref{eq!TsJ}), and the
remainder are more complicated (probably always a sum of two such
products). Mapping this out creates a correspondence from spinor sets to
relations, so from the rows of $M_1$ to its columns; there is obviously
a systematic structure going on here, and nailing it down is an
intriguing puzzle. How this plays out more generally for Kustin--Miller
unprojection \cite{KM}, \cite{PR} and its special cases Tom and Jerry
\cite{TJ} is an interesting challenge.

\subsubsection{Open problems}

To be useful, a structure theory should make some predictions. I hope
that the methods of this paper will eventually be applicable to start
dealing with issues such as the following:
\begin{itemize}
\item $k=3$. A $4\times6$ resolution is a Koszul complex.

\item $k=4$. There are no almost complete intersection Gorenstein
ideals. Equivalently, a $5\times8$ resolution is nonminimal: if $X$ is
Gorenstein codimension~4 and $(L_1,\dots,L_5)$ generate $I_X$ then the
first syzygy matrix $M_1$ has a unit entry, making one of the $L_i$
redundant. This is a well known theorem of Kunz \cite{K}, but I want
to deduce it by my methods.

\item $k=5$. Is it true that a $6\times10$ resolution is a hypersurface
in a $5\times5$ Pfaffian as in \ref{s!exa}, Example~\ref{exa!PfHy}?

The same question for more general odd $k$: are hypersurfaces in a
codimension~3 Gorenstein varieties the only cases? Is this even true for
all the known examples in the literature? This might relate to my even
versus odd remark in \ref{s!evo}.

\item $k=6$. I would like to know whether every case of $7\times12$
resolution is the known Kustin--Miller unprojection from a codimension~4
complete intersection divisor in a codimension~3 complete intersection.

\item $k=8$. As everyone knows, the main case is $9\times16$. How do we
apply the theory to add anything useful to the huge number of known
examples?
\end{itemize}

There are hints that something along these lines may
eventually be possible, but it is not in place yet.

\bigskip

\noindent Miles Reid, \\
Mathematics Institute, University of Warwick, \\
Coventry CV4 7AL, England

\noindent {\it e-mail}: Miles.Reid@warwick.ac.uk

\end{document}